\documentclass[11pt]{amsart}

\newtheorem{thm}{Theorem}[section]
\newtheorem{prop}[thm]{Proposition}
\newtheorem{lem}[thm]{Lemma}

\newtheorem{rem}[thm]{Remark}

\theoremstyle{definition}
\newtheorem{definition}[thm]{Definition}

\theoremstyle{remark}

\numberwithin{equation}{section}

\usepackage{amsmath}
\usepackage{amsthm}
\usepackage{amsfonts}
\usepackage{amssymb}
\usepackage{enumitem}
\usepackage{graphicx,psfrag,subfigure}
\usepackage{epstopdf}
\usepackage{esvect}
\usepackage{epigraph}
\usepackage[spanish, english]{babel}
\usepackage{tikz-cd}

\begin{document}


\title{Fixed points for branched covering maps of the plane}


\author{Alejo Garc\'ia}
\address{Instituto de Matem\'atica y Estad\'istica, Facultad de Ingenier\'ia, Universidad de la Rep\'ublica, 
Montevideo}
\email{alejog@fing.edu.uy}
\urladdr{} 





\begin{abstract}
A well-known result from Brouwer states that any orientation preserving homeomorphism of the plane with no fixed points has an empty non-wandering set. In particular, an invariant compact set implies the existence of a fixed point. In this paper we give sufficient conditions for degree 2 branched covering maps of the plane to have a fixed point, namely:

\textbf{.} A totally invariant compact subset such that it does not separate the critical point from its image
 
\textbf{.} An invariant compact subset with a connected neighbourhood $U$, such that $\mathrm{Fill}(U \cup f(U))$ does not contain the critical point nor its image.

\textbf{.} An invariant continuum such that the critical point and its image belong to the same connected component of its complement. 
\end{abstract}


\maketitle




\section{Introduction}
\label{section:intro}

%






 


\bigskip
The existence of periodic and fixed points for continuous maps of the plane has been extensively studied. A key theorem for the development of this area was given by Brouwer in 1912, \cite{brouwer}:

\begin{thm}\label{thm:brouwer}
Let $f: \mathbb{C} \rightarrow \mathbb{C}$ be an orientation preserving homeomorphism such that $\mathrm{Fix}(f) = \varnothing$. Then every point is wandering. 
\end{thm}

This result triggered a great amount of research, and addresses the simplest case of an open question for plane dynamics: \textit{Does a continuous function of the plane taking a non-separating continuum into itself, necessarily have a fixed point?} \cite{Sternbach}

Cartwright and Littlewood proved in 1951 that \textit{if $f$ is an orientation preseving homeomorphism of the plane with an invariant non-separating continuum, then it has a fixed point in it} (see \cite{Cartwright}). A brief elegant proof of the same result was given two decades later by Brown in \cite{Brown}.

Bell generalized this result to all homeomorphisms in 1978 (see \cite{Bell}), and announced in 1984 that the theorem could be extended to all holomorphic maps (see also \cite{Akis}) -note that these are a particular type \linebreak of \textit{branched covering maps} of the plane, see section \ref{section:preliminares} for definitions-. \linebreak Kuperberg extended the previous theorem for orientation reversing homeomorphisms in 1991, taking out the hypothesis of the continuum being non-separating- \cite{Kuperberg}. 

More recently, the Cartwright-Littlewood theorem was further extended to all orientation preserving branched covering maps of the plane by Fokkink, Mayer, Oversteegen and Tymchatyn in 2007 \cite{Fokkink}.  

We want to find sufficient conditions for branched covering maps of the plane to have fixed (or periodic) points. In particular, we are interested in the next question: Does a branched covering map of the plane with an invariant compact subset, have a fixed point?

In this paper we study the dynamics of degree 2 branched covering maps of the plane. Simple examples can be made for these maps to be periodic point free. See, for example, \cite{Blokh}. 

A similar problem for branched covering maps on the sphere was addressed by Iglesias, Portela, Rovella and Xavier in 2016 (\cite{Iglesias}), proving the next result:

\begin{thm}\label{thm:iprx.esfera}
Let $f$ be a branched covering map of the sphere, with $\lvert d \rvert > 1$. Suppose there exists a simply connected open set $U$, whose closure is disjoint from the set of critical values, and such that $\overline{f^{-1}(U)} \subset U$. Then $f$ has \textit{the rate}.
\end{thm} 

Note that in this case we can restrict the dynamics of $f$ to a branched covering map of the plane -with the same degree-, and get a totally invariant non-separating compact set $K = \displaystyle \bigcap_{n \in \mathbb{N}} f^{-n} (U)$. Note that this set has infinitely many connected components.

By their definition, \textit{the rate} means that $\displaystyle \limsup_{n \rightarrow \infty} \frac{1}{n} log (\#  \mathrm{Fix}(f^n)) \geq log(d).$

Equivalently, it means there is a subsequence of iterates $n_k \rightarrow \infty$ such that $\# \mathrm{Fix}(f^{n_k})$ grows exponentially as $d^{n_k}$, which is the biggest growth to be expected in a non-degenerate case. They actually prove a stronger result, namely, that $f^n$ has at least $d^n$ fixed points, (in particular $\mathrm{Fix}(f) \neq \varnothing$). 

The same mathematicians in 2016 gave sufficient conditions for covering maps of the annulus to have \textit{the rate} (see \cite{Iglesias2}), and as a consequence got the next results, which will be key tools in our proof.

\begin{thm}\label{thm:iprx.anillo}
Let $f$ be a covering map of the annulus, with degree $d > 1$, such that there exists a compact totally invariant subset. Then $f^n$ has at least $d^n - 1$ fixed points. 
\end{thm}

\begin{thm}\label{thm:iprx.esencial}
Let $f$ be a covering map of the annulus, with degree $d > 1$, with an invariant essential continuum. Then $f^{n}$ has at least $d^n - 1$ fixed points in $\mathrm{Fill}(K)$. 
\end{thm}

\section{Notation and preliminaries}
\label{section:preliminares}

Throughout this paper, a \textit{surface} $S$ will be a two dimensional orientable topological manifold. We will say that $S$ is respectively a \textit{plane, annulus} or \textit{sphere} if it is homeomorphic to $\mathbb{C}, \mathbb{C} \backslash \{0\}$ or $\mathbb{S}^2$. A set $U \subset S$ will be a \textit{disc} if it is homeomorphic to $\mathbb{D} = \{z \in \mathbb{C}: {z} < 1\}$. To lighten notation, we will define $m_k: \mathbb{C} \rightarrow \mathbb{C}$ as the map such that $m_k(z) = z^k$. 

All maps considered in this paper are continuous. 

\begin{definition}
Given two points $p, p' \in S$, a \textit{path} or \textit{curve} from $p$ to $p'$ will be a function $\gamma: [0,1] \rightarrow S$ such that $\gamma(0)= p$, $\gamma (1) = p'$. We will say the path is \textit{simple} when $\gamma$ is injective, and \textit{closed} when $p = p'$. 
\end{definition}

We will say $\gamma$ is respectively a \textit{segment, line, circle} if it is homeomorphic to $[0,1], \mathbb{R}, \mathbb{S}^1$. 

\begin{definition}
Given a surface $S$, an \textit{oriented topological foliation} $\mathcal{F}$ is a partition of $S$ in one-dimensional manifolds such that for each $p \in S$, there exists a neighbourhood $U_p$ and a homeomorphism $h : U_p \rightarrow (- 1, 1)\times(- 1, 1)$ which preserves orientation and sends $\mathcal{F}$ into the foliation by vertical lines, oriented from bottom to top. 
\end{definition}

\begin{definition}
Let $U$ be a subset of $\mathbb{C}$, $c \in U$, and take $f: U \rightarrow \mathbb{C}$. We say that $f$ is geometrically conjugate to $m_k$ if there exist foliations $\mathcal{F}$ in $U \backslash \{c\}$ and $\mathcal{F}'$ in $f(U)\backslash \{f(c)\}$, and two homeomorphisms $\phi: U \rightarrow \mathbb{D}$, $\phi':f(U) \rightarrow \mathbb{D}$, such that both $\mathcal{F}$ and $\mathcal{F}'$ are mapped into the radial foliation in $\mathbb{D}$, and such the next diagram commutes:

\begin{center}
	\begin{tikzcd}
    U \arrow{r}{f} \arrow{d}{\phi} & f(U) \arrow{d}{\phi'} \\
     \mathbb{D} \arrow{r}{m_2} & \mathbb{D}
  	\end{tikzcd}
\end{center}
\end{definition}

\begin{definition} 
A \textit{branched covering map} $f: S \rightarrow S$ (or simply \textit{branched covering}), is a map that is a local homeomorphism at each point $p \in S$, except for finitely many \textit{critical points}, each of them having a neighbourhood such that $f$ is geometrically conjugate to $m_k$, with $k \in \mathbb{Z}^+$.
\end{definition}

In that context, the degree of a critical point will be $k$. We will define $\mathrm{Crit}(f)$ as the set of critical points of $f$. Each point which is not the image of a critical point has the same amount of preimages $d$, from which we deduce the \textit{degree} of $f$ is $d$ if $f$ is orientation preserving, and $-d$ if $f$ is orientation reversing. If $\mathrm{Crit}(f) = \varnothing$, we will say $f$ is a \textit{covering map}.

\section{Goals and sketch of the proofs}
\label{section:sketch}

The goal of this paper is to prove the following three results:

\begin{thm}\label{thm:a}
Let $f$ be a degree 2 branched covering map of the plane, and $c \in \mathrm{Crit}(f)$. Suppose there exists a compact set $K$ with $f^{-1}(K) = K$, and such that $c$ and $f(c)$ belong to the same connected component of $\mathbb{C} \backslash K$. Then $\mathrm{Fix}(f) \neq \varnothing$. 
\end{thm}

\begin{thm}\label{thm:b}
Let $f$ be a degree 2 branched covering map of the plane, with a compact invariant set $K$ (not necessarily totally invariant). Suppose there exists $U$ a connected neighbourhood of $K$, such that $\mathrm{Fill}((U) \cup f(U))$ does not contain neither $c$ nor $f(c)$. Then $\mathrm{Fix}(f) \neq \varnothing$. 
\end{thm}

\begin{thm} \label{thm:c}
Let $f$ be a degree 2 branched covering map of the plane, with an invariant continuum $K$, such that it does not separate $c$ from $f(c)$. Then $f$ has a fixed point in $\mathrm{Fill}(K)$. 
\end{thm}

In order to prove these, we will dedicate the next Section to build our main tool: a perturbation $h$ of the map $f$ with \textit{good} properties, given below.

\begin{prop} \label{prop:laproposicion}
Let $f$ be a degree 2 branched covering map of the plane, and $c \in \mathrm{Crit}(f)$. Suppose there exists a compact set $K$ with $f^{-1}(K) = K$, and such that $c$ and $f(c)$ belong to the same connected component of $\mathbb{C} \backslash K$. Suppose in addition that $\mathrm{Fix} (f) = \varnothing$. Then, there exists an orientation preserving homeomorphism $h: \mathbb{C} \rightarrow
 \mathbb{C}$ with $h |_K = \mathrm{Id}$, and $\mathrm{Fix}(h \circ f) = \{ c \} $. 
\end{prop}

\bigskip

Let us outline the roadmap for the proof of this Proposition. 

Given that $c$ is a degree $2$ critical point, we know that there exists a \linebreak neighbourhood of $c$, such that $f$ is locally geometrically conjugate to \linebreak $m_2: \mathbb{C} \rightarrow \mathbb{C}$ -see section \ref{section:preliminares} for definitions-. In Lemma \ref{lem:entorno.maximal} we build a neighbourhood $U$ of $c$ with that property, and such that $\overline{U} \cap \overline{f(U)} \subset \partial U$. We use Lemma \ref{lem:entorno.es.disco} to control the shape of $U$, letting us asume that $U = \mathbb{D}$, $c = 0$. 

In Lemma \ref{lem:dominio.perturbacion} we then construct $V$ the domain of the perturbation $h$ in , which will be a neighbourhood of a path from $c$ to $f(c)$, which is contained in $U \cap f(U)$ except for \textit{one point} (remember that $U \cap f(U) = \varnothing$). We then change coordinates one more time, so that the path becomes a horizontal line -with $c$ on the left of $f(c)$, and the whole neighbourhood $V$ becomes a rectangle, in such a way that we can control the dynamics of $f |_V$. 

The last part of the proof is to define the perturbation in the rectangle -defined as the identity map in its boundary-, which heuristically is the composition of a \textit{vertical contraction} towards the path, composed with a \textit{translation to the left} in each horizontal line.

\section{The good perturbation}
\label{section:theproof}

\begin{lem} \label{lem:entorno.maximal}
There exists a disc $U$ containing $c$, such that

 \textbf{.} $U \cap K = \varnothing$.

    \textbf{.} $f |_U$ is geometrically conjugate to $m_2$.
    
        \textbf{.} $\overline{U} \cap \overline{f(U)} \subset \partial U$. 
\end{lem}

\textit{Proof.} Given that $d(c,K) > 0$, we may take a disc $U_0$ satisfying the first and second property. Let us take a path $\gamma: [0,1] \mapsto \mathbb{R}^2$ from $f(c)$ to $c$, such that $d(\gamma, K) > 0$. Then there exists $V^\varepsilon$ a neighbourhood of $\gamma$, homeomorphic to $[0,1] \times (-\varepsilon,\varepsilon)$ by the homomorphism $h^{\varepsilon}$. We define $V^\varepsilon_t := h^{\varepsilon^{-1}} ([0,t] \times (-\varepsilon,\varepsilon))$, $\tilde{V}_t :=  V^\varepsilon_t \cup f(U_0)$. Let us take $$t_0 = \sup \{t \in [0,1] : \tilde{V}_t \cap  f^{-1} (\tilde{V}_t) = \varnothing\}.$$ Then we may define $U := f^{-1} (\tilde{V}_{t_0})$, and note that the desired intersection is nonempty and is included in the boundary, from which $U$ satifies the third property. \qed

\begin{lem} \label{lem:entorno.es.disco}
Modulo change of coordinates, we may assume that $c = 0$, $U = \mathbb{D}$, and that the radial foliation in $\mathbb{D} \backslash \{0\}$ is mapped 2-1 by $f$ into a foliation by lines in $f(U) \backslash f(0)$. 
\end{lem}

\textit{Proof.} Given that $\partial U$ is a simple closed curve by construction, the Jordan-Schoenflies theorem states a homeomorphism $g: \mathbb{C} \mapsto \mathbb{C}$  such that $g(U) = \mathbb{D}$, $g(\overline{U} ^ c) = \overline {\mathbb{D}} ^ c$ (we may assume $g(c) = 0$). If we define  $f' := g \circ f \circ g^{-1}$, we note it is dynamically conjugate to $f$, from which $f' |_{\mathbb{D}}: \mathbb{D} \rightarrow f'(\mathbb{D})$ is geometrically conjugate to $m_2$, and we obtain a foliation by lines $\mathcal{F}$ in $\mathbb{D} \backslash \{0\}$ as in the end of section \ref{section:preliminares}. Finally, we may conjugate $\mathcal {F}$ to the radial foliation: given that each leaf $\phi \in \mathcal{F}$ gets out of $\mathbb{D}$, simply take the point of the leaf which is in the boundary of $\mathbb{D}$, and define the image of the leaf as the ray which goes through that point. \qed

\bigskip

We will now define the domain of the perturbation, that is, $V \subset \mathbb{R}^2$ such that $h |_{\mathbb{R}^2 \backslash V} = \mathrm{Id}$. We start by taking $U_0 \subset U$ a neighbourhood of $c$ satisfying the first two properties of Lemma \ref{lem:entorno.maximal}

\begin{lem}\label{lem:dominio.perturbacion}
There exists a disc $V$ such that

    \textbf{.} $V \cap K = \varnothing$,
    
    \textbf{.} $V \cap U$ and $V \cap f(U)$ are discs, 
    
    \textbf{.} $U_0 \cup f(U_0) \subset V$,
    
    \textbf{.} $f(V \cap U) \cap V = f(U_0)$,
    
    \textbf{.} $f(V \backslash (U \cup f(U))) \cap V = \varnothing$.
\end{lem}

\textit{Proof.} Let us take $z \in \partial U \cap \partial f(U)$ (given by the lemma \ref{lem:entorno.maximal}). We then take $\hat{z} \in \partial U$ such that $f(\hat{z}) = z$. Let $\gamma$ and $\hat {\gamma}$ be respectively the rays from $z$ to $0$, and from $0$ to $\hat{z}$, parametrized by arc's length. We define $\gamma' := f(\hat{\gamma}). \gamma$. 

For ease, let us call $V^{\varepsilon}$ to the neighbourhoods of $\gamma'$ (as we did with $\gamma$ in Lemma \ref{lem:entorno.maximal}), and define $\gamma'_1 := \gamma' |_{[0,\frac{1}{2}]}$; $\gamma'_2 :=\gamma' |_{[\frac{1}{2},1]}$. Note that there exists $\varepsilon_1$ such that  $V^{\varepsilon_1} \cap K = \varnothing$, because $d(\gamma',K) >0$. 

Given $\mathrm{Fix}(f) = \varnothing$, we know that $f(z) \neq z$. We then get that $z$ is the only point of $\gamma'$ which is outside the discs $U$ and $f(U)$, so there exists $\varepsilon_2$ such that $(V^{\varepsilon_2} \backslash (U \cup f(U))) \cap V^{\varepsilon_2} = \varnothing$.  

Moreover, we have $f(\gamma'_2) \cap \gamma'_1 = f(c)$, from which there also exists $\varepsilon_3$ such that $(f (V^{\varepsilon_3} \cap U) \cap V) \subset f(U_0)$. Finally, taking $\varepsilon_0 := \mathrm{min} \{\varepsilon_1,\varepsilon_2,\varepsilon_3 \}$, we conclude that $V = U_0 \cup V^{\varepsilon_0} \cup f(U_0)$ satisfies the desired properties. \qed

\begin{rem} \label{obs:ubicacion.nuevos.fijos}
$\mathrm{Fix}(h \circ f) \subset \mathrm{Fix}(f) \cup (V \cap f^{-1}(V))$.
\end{rem}

\begin{rem} \label{obs:v.esrectángulo}
With a new change of coordinates, we may assume: 

    \textbf{.} $V$ is the rectangle $(0,10) \times (-1,1)$.

    \textbf{.} $U_0 = (0,2) \times (-1,1)$, $f(U_0) = (8,10) \times (-1,1)$, 
    
    \textbf{.} $c=(1,0)$, $f(c) = (9,0)$.
    
    \textbf{.} $\gamma'$ is the segment between $c$ and $f(c)$.
    
    \textbf{.} $f | _{U_0} (c+(\rho,\theta)) = (f(c) + (f_1 (\rho,\theta),2\theta))$.
    
    \textbf{.} If $(x,y) \in U \cap V$, $0 < x' < x$, then $(x',y) \in U \cap V$. 
 
    \textbf{.} If $(x,y) \in f(U) \cap V$, $x < x' < 10$, then $(x',y) \in f(U) \cap V$. 
\end{rem}

Note that the last two properties in the Remark state that there is a well-defined notion of \textit{left} and \textit{right} inside the rectangle, namely, $U \cap V$ is \textit{on the left} of $f(U) \cap V$. 

\bigskip
We now proceed to build the perturbation $h$. 

\textit{Proof of Proposition \ref{prop:laproposicion}:} We begin by taking $h_1$ a perturbation on the rays of $f(U_0)$, such that $(h_1 \circ f) | _{U_0}$ is, restricted to each ray, the affine map which sends $r_{c,\theta} \cap V$ into $r_{f(c),2\theta} \cap V$. Note that this function may be extended as the identity to the boundary of $f(U_0)$, from which we can extend it to the whole plane, and get $h_1 |{\mathbb{R}^2 \backslash f(U_0)} = \mathrm{Id}$. Observe that the expansion (or contraction) of $h_1 \circ f$ in each ray is uniformly bounded, precisely: $$h_1 \circ f | _{U_0} (c+(\rho,\theta)) = (f(c) + (h_{11}(\theta) \rho,2\theta)) : \sqrt{2}^{-1} \leq h_{11}(\theta) \leq \sqrt{2}$$.

We will now define $h_2$, which is also the identity outside $f(U_0)$, so we will define it on the square $(8,10) \times (-1,1)$. Heuristically, we want a strong contraction towards the horizontal $y=0$, but we need to adjust it so it becomes the identity in the boundary, so we will define a piecewise affine map. We impose simmetry with respect to the line $x=9$, and make it preserve verticals, that is, $h_2 (x,y) = (x, h_{22} (x,y))$, where
\[   
h_{22}(x,y) = 
     \begin{cases}
       \frac {y}{4} &\quad\text{si } \frac{81}{10} \leq x \leq \frac{99}{10} \text{; } \lvert {y} \rvert \leq \frac{1}{2}, \\
       \frac {7y - 3}{4} &\quad\text{si } \frac{81}{10} \leq x \leq \frac{99}{10} \text{; } \lvert {y} \rvert > \frac{1}{2}, \\
       (10x - 99) y + (100 - 10x) \frac{y}{4}  &\quad\text{si } x > \frac{99}{10} \text{; } \lvert {y} \rvert \leq \frac{1}{2}, \\
       (10x - 99) y + (100 - 10x) \frac{7y-3}{4}  &\quad\text{si } x > \frac{99}{10} \text{; } \lvert {y} \rvert > \frac{1}{2}. \\ 
       (81 - 10x) y + (10x - 80) \frac{y}{4} &\quad\text{si } x < \frac{81}{10} \text{; } \lvert {y} \rvert \leq \frac{1}{2}, \\   
       (81 - 10x) y + (10x - 80) \frac{7y-3}{4}  &\quad\text{si } x > \frac{99}{10} \text{; } \lvert {y} \rvert > \frac{1}{2}. \\ 
     \end{cases}
\]
Note that we may extend it as the identity to the boundary of $f(U_0)$.

\bigskip

Finally, we define $h_3$ supported in $V$, preserving horizontals and sending $f(c)$ into $c$ , in a similar fashion to how we defined the map $h_2$. \linebreak Let $h_3 (x,y) = (h_{31} (x,y), y)$, with
\[   
h_{31}(x,y) = 
     \begin{cases}
       \lvert {y} \rvert x + (1-\lvert {y} \rvert) \frac{x}{9} &\quad\text{si } x \leq 9 \\
       \lvert {y} \rvert x + (1-\lvert {y} \rvert) (9x-80) &\quad\text{si } x > 9 \\
     \end{cases}
\]
Note that it may also be extended as the identity yo the boundary of $V$.  \linebreak

Let us define $h := h_3 \circ h_2 \circ h_1$. It only remains to prove:

\begin{lem} \label{lema:nogeneronuevosfijos}
$\mathrm{Fix}(h \circ f) = \{c\}$. 
\end{lem}

\textit{Proof.} Let $w \in \mathrm{Fix}(h \circ f)$. We start by proving $w \in f(U_0)$.

Given that $h$ is supported in $V$, we have that $w \in (V \cap f^{-1}(V))$. By Lemma \ref{lem:dominio.perturbacion}, $f(V \backslash (U \cup f(U))) \cap V = \varnothing$, from which $w \in U$ or $w \in f(U)$.

If $w \in f(U)$, then $f(w) \notin f(U)$, thus we deduce $h \circ f (w) = h_3 \circ f (z)$. Moreover, by the last property of Remark \ref{obs:v.esrectángulo}, and recalling that $h_3$ preserves horizontals, we conclude that $w \notin \mathrm{Fix}(h_3 \circ f)$. On the other hand, if $w \in U $, we use that $f(V \cap U) \cap V = f(U_0)$ (Lemma \ref{lem:dominio.perturbacion}) to conclude $w \in U_0$. 

Let $w = (x,y)$. Note that necessarily $w \in h \circ f (U_0)$. Since $h_1$ and $h_2$ are bijections of $f(U_0)$, we obtain $h\circ f (U_0) = h_3 \circ f (U_0)$, and therefore $w \in (U_0 \cap (h_3 \circ f)(U_0))$, so $\lvert {y} \rvert \leq \frac{\sqrt{2}}{8}$. We then get $w \in U_1 := [0,2] \times [\frac{-\sqrt{2}}{8}, \frac{\sqrt{2}}{8}]$. Let us define $$W_1 := \{ (x,y) \in f(U_0) : x \leq \frac{81}{10}\}; W_2 := \{ (x,y) \in f(U_0) : x \geq \frac{99}{10} \}$$. Observe that: 

\textbf{.} $h_1 \circ f(U_1) \subset [8,10] \times [\frac{-1}{2},\frac{1}{2}]$.

\textbf{.} $f(U_1) \cap W_1 = \varnothing$.

\textbf{.} $h_3 (W_2) \cap U = \varnothing$.

Given that outside the vertical stripes $W_1$, $W_2$, the map $h_2$  divides the height by $4$ in $f(U_1)$, we obtain $h \circ f (U_1) \subset U_2 := [0,2] \times [\frac{-1}{8}, \frac{1}{8}] $, so $w \in U_2$.


Proceeding inductively, we define $U_n := [0,2] \times [\frac{-\sqrt{2}^n}{2^{n+2}}, \frac{\sqrt{2}^n}{2^{n+2}}]$, and conclude $w \in \bigcap_{n \geq 0} U_n$, from which $y=0$. Furthermore, we have that $x \geq 1$ (since $h_1 \circ f (\rho,\theta) = (\rho, 2 \theta)$), so the candidates $w$ are in the interval $I$ of ends \linebreak $c$ and $(2,0)$. There is only left to observe that $f | _I (x,0) = (9x - 8, 0)$ has $c = (1,0)$ as its unique fixed point, which concludes the proof.  \qed

\begin{rem}
The only moment in which we use the fact that $\mathrm{Fix}(f) = \varnothing$ is in Lemma \ref{lem:dominio.perturbacion}, when we ask $z$ not to be fixed. We could exchange it for the (weaker) hypothesis: \textit{``there exists a path from $f(c)$ to $c$ with no fixed point of $f$, contained in $\mathbb{C} \backslash K$''.} The proof works exactly the same way, and in this context we can conclude $\mathrm{Fix} (h \circ f) \backslash \mathrm{Fix} (f) = \{c\}$. As we are searching for fixed points of $f$, both hypothesis are equally as good.
\end{rem}

\section{Proof of the main results} \label{section:applications}

This Section is devoted to the proofs of the three theorems stated in Section \ref{section:sketch}. 

\bigskip

\textit{Proof of Theorem \ref{thm:a}.} The idea is to take a perturbation $h$ as in the Proposition \ref{prop:laproposicion} (without modificating the dynamics in $K$ and such that the critical point becomes fixed, but we do not generate any other fixed point in this process). Given that $h \circ f$ has degree 2, we have that $(h \circ f) ^ {-1} (c) = c$, from which it follows that we can \textit{puncture} the plane in the critical point, and restrict the dynamics to a covering map $g$ of the open annulus $ \mathbb{C} \backslash \{c\}$. 

The set $K$ is then totally invariant for this covering map $g$, so we are on the hypothesis of Theorem \ref{thm:iprx.anillo}, which lets us infer that $\mathrm{Fix}(g) \neq \varnothing$. We then conclude that $\{c\}$ is a proper subset of $\mathrm{Fix}(h \circ f)$, which finishes the proof. \qed

\bigskip

In order to prove Theorem \ref{thm:b}, we will use the following result, which can also be found in \cite{LeCalvez}.

\begin{lem}\label{lem:fijo.en.fill}
Let $f$ be an orientation preserving homeomorphism of the plane, with a compact invariant set $K$. Let $U$ be a connected neighbourhood of $K$. Then $f$ has a fixed point in $V = \mathrm{Fill}(U \cup f(U))$. 
\end{lem}

\textit{Proof.} We know that $f$ has at least one fixed point, by Theorem \ref{thm:brouwer}. Suppose by contradiction that none of them are in $V = \mathrm{Fill}(U \cup f(U))$. Take $S = \mathbb{C} \backslash \mathrm{Fix}(f)$. Given that $\mathrm{Fix}(f)$, is a (closed) invariant set, we may restrict $f$ to $S$. We may also assume that $S$ is connected (otherwise, $f$ preserves the connected component which contains $K$). 

We then take $\tilde{S}$ the universal covering of $S$, which is a plane. Given that $V$ is simply connected, every connected component of $\tilde{\pi}^{-1} (V)$ is homeomorphic to $V$. For the same reason, for any lift $\tilde{f}:\tilde{S} \rightarrow \tilde{S}$ of $f$, and $\tilde{V}$ and any connected component of $\tilde{\pi}^{-1} (V)$, we have that $\tilde{f} (\tilde{V})$ intersects exactly one connected component of $\tilde{\pi}^{-1} (V)$. 

Therefore, if $\tilde{K}$ is the lift of $K$ contained in $\tilde{V}$, we may take a lift $\tilde{f}_0$ of $f$ such that $\tilde{f}_0(\tilde{V})$ only intersects $\tilde{V}$, thus having $\tilde{f}_0(\tilde{K}) = \tilde{K}$. Again by Theorem \ref{thm:brouwer}, we get that $\tilde{f}_0$ has a fixed point, which is a contradiction. \qed

\bigskip

\textit{Proof of Theorem \ref{thm:b}.} We build the same perturbation $h$ as in Proposition \ref{prop:laproposicion}, except in this case we take into account not altering the dynamics in $U$. It is enough to take $U'$ a connected nieghbourhood of $K$ such that it does not contain $c$ nor $f(c)$, and note that there exists $U$ as in the hypothesis, with $(U \cup f(U)) \subset U'$. We then build the perturbation as in Proposition \ref{prop:laproposicion}, such that its domain does not intersect $U'$.

We have that $h \circ f$ has no new fixed points other than $c$.  Let us restrict the dynamics of $h \circ f$ to the annulus $S = \mathbb{C} \backslash \{c\}$, and take the universal covering $\tilde{S}$ (which is a plane). As in Lemma \ref{lem:fijo.en.fill}, any lift $\tilde{V}$ of $V$ is homeomorphic to $V$ (in particular it is bounded), and we can take a lift $\tilde{f}_0$ of $h \circ f$ such that for a lift $\tilde{K}$ of $K$, we get $\tilde{f}_0(\tilde{K})= \tilde{K}$. 

Recalling that $\tilde{f}_0$ is an orientation preserving homemorphism of $\tilde{S}$, we get that it has a fixed point $\tilde{x}$ in $\tilde{V}$, from which we deduce that $x = \tilde{\pi}(\tilde{x})$ is fixed by $h \circ f$ and belongs to $V = \tilde{\pi}(\tilde{V})$, so it is fixed by $f$.  \qed

\begin{rem}\label{rem:nonseparating.sin.c}
If $K$ is a non-separating continuum not containing $c$ nor $f(c)$, the hypothesis of Theorem \ref{thm:b} are immediately verified. Given $K$ is compact and $f$ is continuous, we may control the size of $U \cup \mathrm{Fill}(U)$ and conclude that the fixed point must be found in $K$, thus giving an elementary proof -for degree 2- of a result mentioned in Section \ref{section:intro} and given in \cite{Fokkink}, namely: a positively oriented branched covering map of the plane with an invariant non-separating continuum $K$, has a fixed point in $K$.   
\end{rem}

We now give the proof of the last of our three results.

\bigskip

\textit{Proof of Theorem \ref{thm:c}.} 
We divide the proof in two cases:

\textbf{1.} \textit{The point c belongs to the non-bounded component of $\mathbb{C} \backslash K$.}

We build the perturbation $h$ in the same fashion as in Proposition \ref{prop:laproposicion}. We puncture the plane and go to the annulus $S = \mathbb{C} \backslash c$. As in the last proof, we may take a lift $\tilde{f}_0$ of $h \circ f$ which leaves invariant $\tilde{K}$ a lift of $K$. As in the previous Proof, this implies the existence of a fixed point of $h \circ f$ different from $c$, and so, fixed by $f$. By the same argument used in Remark \ref{rem:nonseparating.sin.c}, we conclude that the fixed point is in $\mathrm{Fill}(K)$.

\textbf{2.} \textit{The point c belongs to a bounded component $C$ of $\mathbb{C} \backslash K$.}

We build the same perturbation, again puncture the plane and take the annulus $A = \mathbb{C} \backslash c$. Observe that $K$ is \textit{essential} in $A$, then by Theorem \ref{thm:iprx.esencial}, we get that $h \circ f |_A$ has a fixed point in $\mathrm{Fill}(K)$, which is also fixed by $f$. Note that in this case we get a little stronger result: the fixed point is actually in $\mathrm{K} \backslash C$ (it may be in $K$ or in one of the bounded components of the complement, in the annulus A).  

\begin{rem}
As in Theorem \ref{thm:b}, we may apply this result when $K$ is an invariant non-separating continuum. In Theorem \ref{thm:c} we actually do not need $K$ to be non-separating, but we do need to impose that $c$ and $f(c)$ belong to the same component of $\mathbb{C} \backslash K$.  
\end{rem}

\section{Room for improvement}

The result proved in this paper leads naturally to other questions, namely:

\begin{itemize}
\item Can we build the perturbation so that it does no alter the set of \textit{periodic points} of $f$? In that case, we will have proved that under the hypothesis of Theorem \ref{thm:a}, $f$ has \textit{the rate}.

\item The same happens if we take Theorem \ref{thm:iprx.anillo} and improve the result to find $f$ has \textit{the rate} in $K$. Keep in mind that the perturbation we built does not modify the dynamics in that set.

\item Can we adjust the techniques to the case where the degree of $f$ is bigger, and get similar results? (This is quite audacious \textit{a priori}, given that the amount of critical points may grow, resulting in cases with more complicated dynamics). 
\end{itemize}


\begin{thebibliography}{Poinc}

\bibitem[Aki99]{Akis}{\ V. Akis.}{\ \textit{On the plane fixed point problem.}}{\ Topology Proc. 24,} (1999), 15-31

\bibitem[Bel78]{Bell}{\ H. Bell.}{\ \textit{A fixed point theorem for plane homeomorphisms.}}{\ Fund. Math. 100 (1978), 119-128.}

\bibitem[BO09]{Blokh}{\ A. M. Blokh and L. G. Oversteegen}{\ \textit{A fixed point theorem for branched covering maps of the plane}}{\ arXiv:0904.2944v1. (2009)}

\bibitem[Brou12]{brouwer}{\ L. E. J. Brouwer.}{\ \textit{Beweis des ebenen Translationssatzes.}}{\ Math. Ann 72} (1912)

\bibitem[Bro77]{Brown}{\ M. Brown.}{\ \textit{A short proof of the Cartwright-Littlewood fixed point Theorem.}}{\  Proc. Amer.Math. Soc. 65, p.372 (1977)}

\bibitem[CL51]{Cartwright}{\ M. L. Cartwright and J. E. Littlewood}{\ \textit{Some fixed point theorems.}}{\ Annals of Math. 54, p.1-37} (1951)

\bibitem[FMOT07]{Fokkink}{\ R. J. Fokkink, J. C. Mayer, L. G. Oversteegen and E. D. Tymchatyn.}{\ \textit{The plane fixed point problem.}}{\ arXiv:0805.1184v2} (2007)

\bibitem[IPRX16]{Iglesias}{\ J. Iglesias, A. Portela, A. Rovella, J. Xavier}{\ \textit{Sphere branched coverings and the growth rate inequality.}}{\ arXiv:1612.02356v1} (2016)

\bibitem[IPRX16.2]{Iglesias2}{\ J. Iglesias, A. Portela, A. Rovella, J. Xavier}{\ \textit{Periodic points of covering maps of the annulus.}}{\ arXiv:1411.5565v2} (2016)

\bibitem[Kup91]{Kuperberg}{\ K. Kuperberg}{\ \textit{Fixed points of orientation reversing homeomorphisms of the plane.}}{\ Proc. of the American Math. Soc. 112}{p.223-229, (1991)}

\bibitem[LeC07]{LeCalvez}{\ P. Le Calvez}{\ \textit{Pourquoi les points périodiques des homéomorphismes du plan tournent-ils autour de certains points fixes?}}{\ Annales scientifiques de l'\'Ecole Normale Sup\'erieure,  Serie 4,  Volume 41}{\ no. 1,  p. 141-176} (2008) 

\bibitem[Ste35]{Sternbach}{\ Sternbach, Problem 107 in}{\ \textit{The scottish Book: MAthematics from the scottish Caf\'e}}{\ Birkhauser, Boston. (1935)}


















\end{thebibliography}
\end{document}